\documentclass[12pt,bezier]{article}
\usepackage{amsmath,amssymb,amsfonts,euscript,graphicx}

\title{   On  rings each of whose finitely generated modules is a direct sum of cyclic modules \thanks
{The research
 of the second author was in part supported by
a grant from IPM (No. 91130413).}
\thanks
 {{\it Key Words}: Cyclic modules; FGC-rings; duo-rings;
 principal ideal rings.}
\thanks {2010{ \it Mathematics Subject Classification}.Primary 16D10, 16D70, 16P20,
  Secondary 16N60. }}

\author{{\small M. Behboodi$^{{\rm a, b}}$\thanks { Corresponding
author;}  and }{\small
 G. Behboodi  Eskandari$^{{\rm a}}$}\\
 {\footnotesize{ $^{\rm a}$Department of Mathematical
 Sciences, Isfahan
University of Technology}}\vspace{-1mm}\\ {\footnotesize{ P.O.Box: 84156-83111, Isfahan, Iran}}\\
 {\footnotesize{ $^{\rm b}$School of Mathematics, Institute for Research in Fundamental Sciences
 (IPM)}}\vspace{-1mm}\\ {\footnotesize{ P.O.Box: 19395-5746, Tehran,
 Iran}}\vspace{-1mm}\\
 {\footnotesize{mbehbood@cc.iut.ac.ir,}}\hspace{3mm}
{\footnotesize{gh.behboodi@math.iut.ac.ir}}}
\date{}
\begin{document}
\maketitle
\begin{abstract}
\noindent{In this paper we study (non-commutative) rings $R$ over
which every finitely generated left module is a direct sum of
cyclic modules (called left FGC-rings). The commutative case was a
well-known problem studied and solved in 1970s by various authors.
It is shown that a Noetherian local left FGC-ring is either an
Artinian principal left ideal ring, or an Artinian principal right
ideal ring, or a prime ring over which every two-sided ideal is
principal as a left and a right ideal. In particular,  it is shown
that a Noetherian local duo-ring $R$ is a left  FGC-ring if and
only if $R$ is a right FGC-ring, if and only if,  $R$ is a
principal ideal ring. Moreover,   we obtain that if $R=\Pi_{i=1}^n
R_i$ is a finite product of Noetherian  duo-rings $R_i$ where each
$R_i$ is prime or local, then $R$ is a left FGC-ring if and only
if $R$ is a principal ideal ring. }
\end{abstract}

{\bf 1. Introduction}

 \noindent The question of which commutative rings have the property that
every finitely generated module is a direct sum of cyclic modules
has been around for many years. We will call these rings
FGC-rings. The problem originated in I. Kaplansky's papers
\cite{Kap1} and \cite{Kap2}, in which it was shown that a local
domain is FGC if and only if it is an almost maximal valuation
ring. For several years, this is one of the major open problems in
the theory. R. S. Pierce \cite{Pie} showed that the only
commutative FGC-rings among the commutative (von Neumann) regular
rings are the finite products of fields. A deep and difficult
study was made by Brandal \cite{Bra2}, Shores-R. Wiegand
\cite{Sho}, S. Wiegand \cite{Wie2}, Brandal-R. Wiegand \cite{Bra3}
and V$\acute{a}$mos \cite{Vam}, leading to a complete solution of
the problem in the commutative case. To show that a commutative
FGC-ring cannot have an infinite number of minimal prime ideals
required the study of topological properties (so-called Zariski
and patch topologies). For complete and more leisurely treatment
of this subject, see Brandal \cite{Bra1}. It gives a clear and
detailed exposition for the reader wanting to study the subject.
The main result reads as follows: A commutative ring $R$ is an
FGC-ring exactly if it is a finite direct sum of commutative rings
of the following kinds: (i) maximal valuation rings; (ii) almost
maximal B$\acute{e}$zout domains; (iii) so-called torch rings (see
\cite{Bra1} or \cite{Fac} for more details on the torch rings).

The corresponding problem in the non-commutative case is  still
open; see \cite[Appendix B. Dniester Notebook: Unsolved Problems
in the Theory of Rings and Modules. Pages 461-516]{Sab}   in which the following problem is considered.\vspace{3mm}\\
 {\bf \cite[Problem 2.45]{Sab}}  (I. Kaplansky, reported by A. A.
 Tuganbaev):
{\it Describe the rings in which every one-sided ideal is
two-sided and over which every finitely generated module can be
decomposed as a direct sum of cyclic modules.}

Through this paper,  all rings have identity elements and all
modules are unital. A {\it left FGC-ring}  is a ring $R$ such that
each finitely generated  left $R$-module is a direct sum of cyclic
submodules. A {\it right  FGC-ring}  is defined similarly, by
replacing the word left with right above.  A ring $R$ is called a
{\it FGC-ring} if it is a both left and right FGC-ring.  Also, a
ring $R$ is  called {\it duo-ring} if each one-sided ideal of $R$
is two-sided. Therefore, the Kaplansky problem is: {\it Describe
the FGC-duo-rings.}

In this paper, we study left FGC-rings and, among other results,
we will present a partial solution to the above problem of
Kaplansky.\\

\noindent{\bf 2. On Left FGC-Rings}

\noindent A ring $R$ is {\it local} in case $R$ has a unique left
maximal ideal. An Artinian (resp. Noetherian) ring is a ring which
is both left and right Artinian (resp. Noetherian). A {\it
 principal ideal ring}  is a ring which is both left and right principal ideal
 ring. Also, for a subset $S$ of $_RM$, we denote by l.Ann$_R(S)$,
the left annihilator of $S$ in $R$. A left $R$-module $M$ which
has a composition series is called a module of finite length. The
length of a composition series of $_RM$ is said to be the length
of $_RM$ and denoted by length$(_RM)$.

Note that in the sequel, if we have proved certain results for
rings or modules ``on the left," then we shall use such results
freely also ``on the right," provided that these results can
indeed be proved by the same arguments applied ``to the other
side."

 We begin with the following lemma which is an associative, non-commutative version of Brandal [2, Proposition 4.3]
  for local rings $(R, {\mathcal{M}})$  with
  ${\mathcal{M}}^{2}=(0)$. Also, the proof is based on a slight
  modification of the proof of \cite [Theorem 3.1]{Beh}.

\noindent{\bf  Lemma  2.1.} {\it Let $(R, {\mathcal{M}})$
 be a  local ring with  ${\mathcal{M}}^{2}=(0)$ and $_{R}{\mathcal{M}}=Ry_{1}\oplus\ldots\oplus Ry_{t}$
such that $t\geq 2$ and each $Ry_i$ is a minimal left ideal of
$R$. If there exist $0\neq x_1$, $x_2\in {\mathcal{M}}$ such that
$x_1R\cap x_2R=(0)$, then the left $R$-module $(R\oplus R)/R(x_1,
x_2)$ is not a direct sum of cyclic modules.}

\noindent{\bf Proof.} Since
$_{R}{\mathcal{M}}=Ry_{1}\oplus\ldots\oplus Ry_{t}$ and
 each $Ry_i$ is a minimal left ideal of $R$, we conclude that  $R$ is of finite composition length and
length$(_{R}R)=t+1$. We put $_{R}G=(R\oplus R)/R(x_1, x_2)$. Since
$x_1$, $x_2\in {\mathcal{M}}$ and ${\mathcal{M}}^{2}=(0)$, we
conclude that ${\rm l.Ann}_R(R(x_1, x_2))={\mathcal{M}}$. Thus
$R(x_1, x_2)$ is simple and hence

   $~~~~~{\rm length}(_{R}G)= 2 \times
{\rm length}(_{R}R)-{\rm length}(_{R}R(x_1, x_2))=2(t+1)-1.$

 We claim that every non-zero cyclic submodule $Rz$ of $G$ has length $1$ or $t+1$.
If ${\mathcal{M}}z=0$, then length$(Rz)=1$ since $Rz\simeq
R/{\mathcal{M}}$. Suppose that ${\mathcal{M}}z\neq 0$, then there
exist $c_1,c_2\in R$ such that $z=(c_1, c_2)+R(x_1, x_2)$. If
$c_1,c_2\in {\mathcal{M}}$, then ${\mathcal{M}}z=0$, since
${\mathcal{M}}^2=0$. Thus without loss of generality,  we can
assume that $z=(1, c_2)+ R(x_1, x_2)$ (since if $c_1\not\in
{\mathcal{M}}$, then $c_1$ is unit). Now let $r\in {\rm
l.Ann}_R(z)$, then $r(1, c_2)=t(x_1, x_2)$ for some $t\in R$. It
follows that $r=tx_1$ and $rc_2=tx_2$. Thus $tx_2=tx_1c_2$. If
$t\notin {\mathcal{M}}$, then $t$ is unit and so $x_2=x_1c_2$ that
it is contradiction (since $x_1R\cap x_2R=(0)$). Thus $t\in
{\mathcal{M}}$ and so $r=tx_1=0$. Therefore, ${\rm l.Ann}_R(z)=0$
and so $Rz\cong R$. It follows that length$(Rz)=t+1$.

 Now suppose the assertion of the lemma is false. Then $_RG$  is a direct sum of
cyclic modules and since  $_RG$ is of finite length,  we have

$~~~~~~~~~~~~~~~~~~~~~G=Rw_1\oplus\ldots\oplus Rw_k\oplus
Rv_1\oplus \ldots\oplus Rv_l,$

\noindent where $l, k\geq 0$, and each  $Rw_i$ is of length $t+1$
and each $Rv_j$ is of length $1$. Clearly ${\mathcal{M}}\oplus
{\mathcal{M}}$ is not a simple  left $R$-module.
 Since $R(x_1, x_2)$ is simple,
  ${\mathcal{M}}G=({\mathcal{M}}\oplus {\mathcal{M}})/R(x_1,x_2)\neq 0$. It follows
 that $k\geq 1$.  Also, length$(_{R}G)=2(t+1)-1=k(t+1)+l$ and this implies that $k=1$ and $l=t$. Since ${\mathcal{M}}v_i=0$ for each $i$, ${\mathcal{M}}G={\mathcal{M}}w_1$
    and hence

$~~~~~~~~~~~~~~~~~~~~~G/{\mathcal{M}}G\simeq
Rw_1/{\mathcal{M}}w_1\oplus Rv_1\oplus \ldots\oplus Rv_t$.

 It follows that ${\rm length}(_{R}G/{\mathcal{M}}G)=1+t$.
On the other hand, we have $$G/{\mathcal{M}}G\cong
R/{\mathcal{M}}\oplus R/{\mathcal{M}}$$ and so
length$(_{R}G/{\mathcal{M}}G)=2$ and so $t=1$,  a contradiction.
Thus the left $R$-module $(R\oplus R)/R(x_1,
x_2)$ is not a direct sum of cyclic modules.$~\square$\\

We recall that the socle ${\rm soc}(_RM)$ of a left module $M$
over a ring $R$ is defined to be the sum of all simple submodules
of $M$.

\noindent{\bf  Theorem  2.4.} {\it Let $(R, {\mathcal{M}})$
 be a  local ring such that $_R{\mathcal{M}}$ and  ${\mathcal{M}}_R$ are finitely generated.
  If  every left $R$-module with two generators is a direct sum of
cyclic modules,    then  either ${\mathcal{M}}$  is a principal
left ideal  or   ${\mathcal{M}}$  is a  principal right ideal.}

\noindent{\bf Proof.}  We can assume that ${\mathcal{M}}$ is not a
principal left ideal of $R$. One can easily see that
${\mathcal{M}}_R$  is generated by $\{x_1,\cdots, x_n\}$ if and
only if
   ${\mathcal{M}}/{\mathcal{M}}^2$ is generated by the set
 $\{x_1+{\mathcal{M}}^2,\cdots, x_n+{\mathcal{M}}^2\}$
 as a  right  ideal of  $R/{\mathcal{M}}^2$. Thus    it suffices
to show that ${\mathcal{M}}/{\mathcal{M}}^{2}$ is a principal
right ideal of  $R/{\mathcal{M}}^{2}$. Since every left $R$-module
with two generators is a direct sum of cyclic modules, we conclude
that every left $R/{\mathcal{M}}^{2}$-module with two generators
is also a direct sum of cyclic modules. Therefore, without loss of
generality we can assume that ${\mathcal{M}}^{2}=(0)$.  It follows
that ${\rm soc}(_RR)={\rm soc}(R_R)={\mathcal{M}}$.   Since
$_R{\mathcal{M}}$ is finitely generated,
$_{R}{\mathcal{M}}=Ry_{1}\oplus\ldots\oplus Ry_{t}$ such that
$t\geq 2$ and each $Ry_i$ is a minimal left ideal of $R$. We claim
that ${\mathcal{M}}_R=xR$, for if not, then we can assume that
${\mathcal{M}}_{R}=\oplus_{i\in I}x_{i}R$ where $|I|\geq 2$ and
each $x_{i}R$ is a minimal right ideal of  $R$. We can assume that
$\{1,2\}\subseteq I$ and so $0\neq x_1$, $x_2\in {\mathcal{M}}$
and  $x_1R\cap x_2R=(0)$. Now by Lemma 2.1, the left $R$-module
$(R\oplus R)/R(x_1, x_2)$ is not a direct sum of cyclic modules, a
contradiction. Thus ${\mathcal{M}}$ is principal  as a right ideal
of $R$.$~\square$

A  ring whose lattice of left ideals is linearly ordered under
inclusion, is called a {\it left uniserial} ring. A  {\it
uniserial ring}
  is a ring which is both left and right uniserial. Note that left and right uniserial rings are
  in particular local rings and commutative uniserial rings are also known as valuation rings.

Next, we need the following lemma from \cite{Nic}.

\noindent{\bf  Lemma  2.5.}  (See Nicholson and S\'{a}nchez-Campos
\cite[Theorem 9]{Nic}) {\it  For any
ring $R$, the following statements  are equivalent:}\vspace{2mm}\\
(1)  {\it $R$ is local,  $J(R)=Rx$ for some $x\in R$ and $x^k=0$ for some  $k\in\Bbb{N}$.}\vspace{1mm}\\
  (2) {\it There exist $x\in R$ and $k\in\Bbb{N}$ such that  $x^{k-1}\neq 0$ and
  $R\supset Rx\supset \ldots\supset Rx^k=(0)$ \indent  are the only left ideals of $R$.}\vspace{1mm}\\
   (3) {\it $R$  is left uniserial of finite composition
   length.}\vspace{3mm}

\noindent{\bf  Theorem  2.6.} {\it Let $(R, {\mathcal{M}})$
 be a local ring  such that $_R{\mathcal{M}}$ and  ${\mathcal{M}}_R$ are finitely generated and ${\mathcal{M}}^k=(0)$
 for some $k\in\Bbb{N}$.
  If  every left $R$-module with two generators is a direct sum of
cyclic modules,    then  either  $R$ is a left Artinian principal
left ideal ring   or   $R$ is a right  Artinian  principal right
ideal ring.}

\noindent{\bf Proof.} Assume that every left $R$-module with two
generators is a direct sum of cyclic modules. Then by Theorem 2.4,
either ${\mathcal{M}}$ is a principal left ideal  or
${\mathcal{M}}$ is a  principal right ideal. If ${\mathcal{M}}$ is
a principal left ideal, then by Lemma  2.5,  $R$ is a left
Artinian principal left ideal ring. Thus we can assume that
${\mathcal{M}}$ is a principal right ideal.   Then by using  Lemma
2.5 to the right side, $R$ is a right Artinian principal right
ideal ring.$~\square$

Next, we need the following lemma from Mohamed H. Fahmy-Susan
Fahmy\cite{Fah}. We note that their definition of a local ring is
slightly different than ours; they defined a {\it local ring}
(resp. {\it scalar local ring}) as a ring $R$  such that it
contains a unique maximal ideal ${\mathcal{M}}$ and
$R/{\mathcal{M}}$ is an Artinian ring (resp. division ring). Thus
our definition of a {\it local ring} and  the  {\it scalar local
ring} coincide.

\noindent{\bf  Lemma  2.7.} (See \cite[Theorem 3.2]{Fah} {\it  Let
$(R, {\mathcal{M}})$ be non-Artinian Noetherian local ring. Then
the following conditions are equivalent:}\vspace{2mm}\\
(1) {\it ${\mathcal{M}}$ is principal as a right ideal.}\\
 (2) {\it ${\mathcal{M}}$ is principal as a left ideal.}\\
  (3) {\it Every two-sided ideal of $R$ is  principal as a left ideal.}\\
   (4) {\it Every two-sided ideal of $R$ is  principal as a right ideal.\vspace{2mm}\\
    \indent Moreover, $R$ is a prime ring.}

Now we are in a position to prove the main theorem of this
section.

\noindent{\bf  Theorem  2.8.} {\it Let $(R, {\mathcal{M}})$
 be a Noetherian local ring.  If every left $R$-module with two generators is a direct sum of
cyclic modules,  then
  one of the following holds:}\vspace{2mm}\\
 (a) {\it $R$ is an Artinian  principal left ideal ring.}\\
 (b) {\it $R$ is an   Artinian  principal right ideal ring.}\\
 (c) {\it $R$ is a prime ring and every two-sided ideal of $R$ is  principal as both left and \indent right
 ideals.}

\noindent{\bf Proof.}  First we assume that $R$ is  an Artinian
ring. Thus by Theorem 2.6, either $R$ is an Artinian  principal
left ideal ring or  $R$ is an   Artinian  principal right ideal
ring. Now we assume that $R$ is not an Artinian ring. By Theorem
2.4, either  ${\mathcal{M}}$  is a principal left ideal  or
${\mathcal{M}}$  is a  principal right ideal. Thus  by Lemma 2.7,
$R$ is a prime ring and every two-sided ideal of $R$ is principal
as both left and right  ideals.$~\square$\\

\noindent{\large\bf 4. A Partial Solution of  Kaplansky's Problem
on Duo-Rings}

\noindent A ring $R$ is said to be {\it left } (resp. {\it right})
{\it hereditary}  if every left  (resp. right) ideal of $R$ is
projective as a left  (resp. right) $R$-module. If $R$ is both
left and right hereditary, we say that $R$ is hereditary. Recall
that a PID  is a domain $R$ in which any left and any right ideal
of $R$ is principal.  Clearly, any PID  is hereditary.

Let $R$ be an hereditary prime ring with quotient ring $Q$ and $A$
be a left  $R$-module. Following Levy \cite{Lev}, we say that
$a\in A$  is a {\it torsion element}  if there is a regular
element $r\in  R$ such that $ra=0$. Since, by Goldie's theorem,
$R$ satisfies the Ore condition, the set of torsion elements of
$A$ is a submodule $t(A)\subseteq  A$. $A/t(A)$ is evidently
torsion free (has no torsion elements).

\noindent{\bf  Lemma 3.1.} (Eisenbud-Robson \cite[Theorem
2.1]{Eis1}) {\it Let $R$ be an hereditary Noetherian prime ring,
and let $A$ be a finitely generated left $R$-module. Then $A/t(A)$
is projective and $A\cong t(A)\oplus A/t(A)$.}

A {\it Dedekind prime  ring} \cite{Rob} is an hereditary
Noetherian prime ring with no proper idempotent two-sided ideals
(see \cite{Eis2}). Clearly if a duo-ring $R$ is a  PID, then $R$
is a Dedekind prime ring.

\noindent{\bf  Lemma 3.2.} (Eisenbud-Robson \cite[Theorem
3.11]{Eis1}) {\it  Let $R$ be a  Dedekind prime ring. Then every
finitely generated torsion left $R$-module $A$ is a direct sum of
cyclic modules.}

\noindent{\bf  Lemma 3.3.} (Eisenbud-Robson \cite[Theorem
2.4]{Eis1}) {\it  Let $R$ be a Dedekind prime ring, and let $A$ be
a
projective left $R$-module. Then:}\vspace{2mm}\\
 (i) {\it If $A$ is finitely generated, then $A\cong F\oplus I$ where $F$ is a finitely generated free
 \indent  module  and $I$ is a
 left ideal  of $R$.}\\
  (ii) {\it If $A$ is
not finitely generated, then $A$ is free.}

\noindent{\bf Proposition 3.4.} {\it Let $R$ be a Dedekind prime
ring. If $R$ is a left principal ideal ring, then $R$ is a left
FGC-ring.}

\noindent{\bf Proof.} Suppose that $A$ is a finitely generated
left $R$-module. Since $R$ is a Dedekind prime ring,  $R$ is
Noetherian and so $A$ is also a Noetherian left $R$-module. Thus
by Lemma 3.1, $A/t(A)$ is projective and $A\cong t(A)\oplus
A/t(A)$. By Lemma 3.2, $t(A)$ is a direct sum of cyclic modules.
Also by Lemma 3.3, $A/t(A)\cong F\oplus I$ where $F$ is a  free
module and $I$ is a  left ideal  of $R$. Since $R$ is a  principal
left ideal ring, $I$ is a cyclic left $R$-module, i.e.,  $A/t(A)$
is a direct sum of cyclic modules. Thus, $A\cong t(A)\oplus
A/t(A)$ is a direct sum of cyclic modules. Therefore, $R$ is a
left FGC-ring.$~\square$

The following proposition  is an answer to the question: ``What is
the class of FGC Noetherian prime  duo-rings?"

\noindent{\bf Proposition  3.5.} (See also Jacobson \cite[Page 44,
Theorems 18 and 19]{Jac}) {\it Let $R$ be a Noetherian prime
duo-ring (i.e., $R$ is a Noetherian duo-domain). Then the
following statements are equivalents:}\vspace{2mm}\\
(1) {\it $R$ is an FGC-ring.}\\
  (2) {\it $R$ is a left  FGC-ring.}\\
 (3) {\it  $R$ is a  principal ideal ring.}\vspace{2mm}\\
 {\it The same characterizations also apply for right $R$-modules.}

\noindent{\bf Proof.} (1) $\Rightarrow$ (2) is clear.\\
(2) $\Rightarrow$ (3).  Suppose that $I$ is an ideal of $R$. Since
$I$ is a direct sum of principal ideals of $R$
and $R$ is a domain, we conclude that $I$ is principal. Thus, $R$ is a principle ideal ring.\\
(3) $\Rightarrow$ (1) is by Proposition 3.4.$~\square$

A {\it left} (resp., {\it right})  {\it K\"{o}the ring}  is a ring
$R$ such that each left (resp.,  right) $R$-module is a direct sum
of cyclic submodules.   A ring $R$ is called a {\it K\"{o}the
ring} if it is a both left and right K\"{o}the ring. In \cite{Kot}
K\"{o}the proved that  an Artinian principal ideal ring is a
K\"{o}the ring. Furthermore, a commutative  ring $R$  is a
K\"{o}the ring if and only if $R$ is an Artinian principal ideal
ring (see Cohen and Kaplansky \cite{Coh}). The corresponding
problem in the non-commutative case is  still open (see
\cite[Appendix B, Problem 2.48]{Sab} or  Jain-Srivastava
\cite[Page 40, Problem 1]{Jai}. Recently,  a generalization of the
K\"{o}the-Cohen-Kaplansky theorem is given in   \cite{Beh}. In
fact: in  \cite[Corollary 3.3.]{Beh}, it is shown that if  $R$ is
a ring in which all idempotents are central, then $R$ is a
K\"{o}the ring if and only if   $R$ is an Artinian principal ideal
ring.

 Next, the following theorem is an answer to the question:
``What is the class of FGC Noetherian local duo-rings?"

\noindent{\bf  Theorem 3.6.} {\it Let $(R, {\mathcal{M}})$
 be a  Noetherian local duo-ring. Then the following statements are equivalent:}\vspace{2mm}\\
 (1) {\it $R$ is an FGC-ring.}\\
  (2) {\it $R$ is a left  FGC-ring.}\\
 (3) {\it Every left $R$-module with two generators is a direct sum of
cyclic modules.}\\
 (4) {\it Either $R$ is an Artinian  principal ideal ring or $R$ is a principal ideal
 domain.}\\
 (5) {\it  $R$ is a  principal ideal ring.}\vspace{2mm}\\
 {\it The same characterizations also apply for right $R$-modules.}

   \noindent{\bf Proof.} (1) $\Rightarrow$ (2)  $\Rightarrow$ (3) is clear.\\
   (3) $\Rightarrow$ (4).   Suppose that every left $R$-module with two generators is a direct sum of
cyclic modules. Thus by Theorem 2.4,   ${\mathcal{M}}$  is
principal as both left and right  ideals. If $R$ is Artinian, then
by Theorem 2.6, $R$ is an Artinian principal ideal ring. If $R$ is
not Artinian, then by
Lemma 2.7, $R$ is a principal ideal domain. \\
(4) $\Rightarrow$ (1). If $R$ is  an Artinian principal ideal
ring, then by the K\"{o}the result,  each left, and each right
$R$-module is a direct sum of cyclic modules. Thus  $R$ is an
FGC-ring. Now
assume that $R$ is a principal ideal domain. Then by Proposition 3.5, $R$ is an FGC-ring.\\
(4) $\Rightarrow$ (5)  is clear.\\
(5) $\Rightarrow$ (4). Assume that  $R$ is a  principal ideal
ring. Then ${\mathcal{M}}$  is  principal as both left and right
ideals. If $R$ is Artinian, then Lemma 2.5,  $R$ is an Artinian
principal ideal ring. If $R$ is not Artinian, then by Lemma 2.7,
$R$ is a principal ideal domain.$~\square$

Let $R=\Pi_{i=1}^n R_i$ be  a finite product of rings $R_i$.
Clearly  $R$ is a principal ideal ring if and only if each $R_i$
is a principal ideal ring. On the other hand if $R$ is a left
FGC-ring, then each $R_i$ is also a left FGC-ring. Thus as a
corollary of Proposition 3.5 and  3.6, we have the following
result.

\noindent{\bf  Corollary   3.7.} {\it  Let $R=\Pi_{i=1}^n R_i$
 be a finite product of Noetherian  duo-rings $R_i$ such that each $R_i$ is a domain or a local ring.  Then the  following statements are equivalent:}\vspace{2mm}\\
(1) {\it $R$ is an FGC-ring.}\\
  (2) {\it $R$ is a left  FGC-ring.}\\
 (3)  {\it  $R$ is a  principal ideal ring.}\vspace{2mm}\\
 {\it The same characterizations also apply for right $R$-modules.}

Next, we need the following   lemma from \cite{Hab} about Artinian
duo-rings (its proof is worthwhile even in the commutative case
(see \cite[Corollary 4]{Hab} or \cite[Lemma 4.2]{Kar})

\noindent{\bf Lemma 3.8.}  {\it Let $R$ be an Artinian duo-ring.
Then $R$ is a finite direct product of Artinian local duo rings.}

Next, we give the  following characterizations of an Artinian FGC
duo-ring. In fact, on Artinian duo-rings, the notions ``FGC" and
``K\"{o}the" coincide.

\noindent{\bf  Theorem 3.9.} {\it Let $R$
 be an Artinian Duo-ring. Then the following statements are equivalent:}\vspace{2mm}\\
 (1) {\it $R$ is a left FGC-ring}\\
 (2) {\it $R$ is an  FGC-ring}\\
 (3) {\it Every left $R$-module with two generators is a direct sum of
cyclic modules.}\\
 (4) {\it $R$ is a left  K\"{o}the-ring}\\
 (5) {\it $R$ is a K\"{o}the-ring}\\
 (6) {\it $R$ is a  principal ideal ring.}\vspace{2mm}\\
 {\it The same characterizations also apply for right $R$-modules.}

\noindent{\bf Proof.} Since $R$ is an Artinian duo-ring, by Lemma
3.8,  $R=\Pi_{i=1}^n R_i$ such that each $R_i$ is  an Artinian
local duo-ring.   Thus by the     K\"{o}the result  and Corollary
3.7, the proof is complete.$~\square$

\end{document}